\documentclass[12pt,a4paper]{article}
\usepackage{latexsym,amssymb}
\usepackage{t1enc}

\usepackage{amsmath}
\usepackage{amssymb}
\usepackage{amsthm}
\usepackage[cp1250]{inputenc}
\usepackage[T1]{fontenc}
\usepackage[english]{babel}
\usepackage[dvips]{graphicx}
\usepackage[a4paper, left=3.2cm, right=2.8cm, top=2.3cm, bottom=3.4cm, headsep=1.2cm]{geometry}

\newtheorem{theorem}{Theorem}[section]
\newtheorem{lemma}[theorem]{Lemma}
\newtheorem{proposition}[theorem]{Proposition}
\newtheorem{corollary}[theorem]{Corollary}
\newtheorem{remark}[theorem]{Remark}

\def\R{\mathbb{R}}
\def\C{\mathbb{C}}
\def\N{\mathbb{N}}

\def\d{\mathrm{\,d}}
\def\lma{\lambda}
\def\part{\partial}

\begin{document}

\title{\textbf{Continuity and averaging for parabolic evolution systems}}
\author{ Aleksander \'{C}wiszewski\footnote{The research was supported by the NCN Grant 2013/09/B/ST1/01963.} , Renata {\L}ukasiak\\
\emph{Nicolaus Copernicus University}\\
\emph{ul. Chopina 12/18, 87-100 Toruń, Poland}\\
e-mail: Aleksander.Cwiszewski@umk.pl, Renata.Lukasiak@umk.pl }

\maketitle

\begin{abstract}
Averaging principle for abstract non-autonomous parabolic evolution equations governed by time-dependent family of positive sectorial operators is proved. Apart from linear case also a nonlinear version for continuous perturbations is provided. These results extend Henry's averaging technique \cite{Henry} to the case where linear part of the equation depends on time and results due to Ilyin \cite{Ilyin}.
\end{abstract}

\indent \indent {\bf Key words}: evolution system, parabolic equation, averaging

\indent \indent {\bf MSC:} 47D06, 34G05, 34G20

\section{Introduction}\
\indent We are interested in the limit behavior when $\lambda \rightarrow 0^+$ of solutions of the following family of equations
\begin{equation}
\tag{$P_\lambda$}\label{eq:introduction}
\dot{u}(t) = - A(t /\lma) u(t)+ F(t/ \lma, u(t)), \ t>0,
\end{equation}
\noindent
where, for each $t \geq 0$, $A(t)$ is a sectorial operator in a Banach space $X$,   $F: [0,+\infty) \times X^\alpha \rightarrow X$ is a continuous mapping with sublinear growth and $\lambda>0$ is a parameter.
Here $X^\alpha$, $\alpha \in [0,1)$, is the fractional power space determined by the operators $A(t)$, $t\geq 0$ (under our assumptions it is independent of $t$).\\
\indent The classical Bogolyubov-Krylov-Mitropol'skii averaging idea was applied to ordinary differential equations. Generally speaking, it states that solutions of $\dot u (t)=f(t/\lma, u(t))$ converge as $\lma\to 0^+$ to solutions of the averaged equation $\dot u(t) = \widehat f(u(t))$ where $\widehat f$ is the average of $f$, i.e.
$\widehat f(u)=\lim_{\omega\to +\infty} \frac{1}{\omega}\int_{0}^{\omega} f(t,u)\d t$, which exists under proper assumptions (see \cite{Bogoliubov}).  In other words, the averaging scheme enables to study the asymptotic behavior of solutions of nonautonomous problems in terms of the related averaged one.\\ 
\indent The aim of this paper is to extend the averaging principle to abstract infinite dimensional evolution equations with unbounded linear part depending on time and apply it to parabolic partial differential equations with elliptic differential operator depending on time.  One of first general averaging principles in infinite dimension was derived by Henry in \cite{Henry} for problems
$$
\dot u = - Au +F(t/\lma, u),
$$
where $-A$ is a sectorial operator, i.e. for \eqref{eq:introduction}
with $A(t)\equiv A$ (the linear part is independent of time). These results were applied to parabolic partial differential equations.  Equations governed by perturbations of $C_0$-semigroup generators  where studied by Hale and Verduyn Lunel in \cite{Hale-Lunel} and applied their results to retarded differential equations. The existence and structure of global attractor for dissipative systems were considered in \cite{Ilyin} where examples for reaction-diffusion, a two-dimensional Navier-Stokes system and a damped wave equation may be found. Averaging principle for the case where $-A$ is the generator of a $C_0$-semigroup of strict contractions was considered by Couchouron and Kamenskii in \cite{Couchouron-Kamenskii}
where $F$ was an upper semicontinuous set-valued map being a $k$-set contraction with respect to a measure of noncompactness. Antocci and Prizzi in \cite{Antoci-Prizzi} studied averaging and the upper semicontinuity of attractors of parabolic partial parabolic equations on $\R^N$ with elliptic operator having time dependent coefficients.  Prizzi in \cite{Prizzi} considers invariant sets of parabolic equations on unbounded domains. An improvement of Henry's general result was given in \cite{Cwiszewski-RL} where initial values could converge in $X$ (not necessarily in $X^\alpha$). Averaging appeared useful in studying periodic solutions -- see \cite{Kamenskii}, \cite{Cwiszewski-JDE-2005}, \cite{Cwiszewski-Kokocki-DCSA2008}, \cite{Cwiszewski-Kokocki-JEE2010}, \cite{Cwiszewski-CEJM}, \cite{Kokocki1}, \cite{Kokocki2} as well as \cite{Cw-Luk-2015}.\\
\indent Here we provide a general  approach for abstract evolution equations of parabolic type that embrace the case considered in \cite{Antoci-Prizzi} and extend Henry's general theory.  Our intention was to come up with results that have minimal assumptions and can be verified in practice for partial differential equations with time dependent coefficients. To this end we adapt the so-called {\em parabolic evolution systems} obtained by use of the Levy construction (see \cite{Pazy} and \cite{Tanabe}). We only add an assumption that allows us to consider fractional space $X^\alpha$ that is independent of $t$ and some specific conditions that are necessary for averaging.\\
\indent Throughout the whole paper the family $\{ A(t)\}_{t\geq 0}$ of operators in a Banach space $X$ is assumed to satisfy the following conditions:
\begin{enumerate}
\item [$(\mathcal{P}_0)$]
for each $t\geq 0$, the operator $A(t)$ is sectorial
and for some fixed $\alpha \in (0,1)$ the associated fractional spaces $X^\alpha$ are independent of $t$ and have equivalent norms, i.e there are $C_{1,\alpha}, C_{2,\alpha}>0$ such that, for any $t \geq 0$
$$
C_{1,\alpha}\|u\|_\alpha \leq \|A(t)^\alpha u\| \leq C_{2,\alpha}\|u\|_\alpha, \textnormal{ for all } u \in X^\alpha,
$$
where $A(t)^\alpha$ is the fractional power of $A(t)$ and $\|u\|_\alpha:= \|A(0)^\alpha u\|$;
\item [$(\mathcal{P}_1)$]
for each $t\geq 0$ the domain $X^1$ of $A(t)$ is independent of $t$ and there are $C_{1}, C_{2}>0$ such that, for any $t \geq 0$
$$
C_{1}\|u\|_1 \leq \|A(t)u\| \leq C_{2}\|u\|_1, \textnormal{ for all } u \in X^1,
$$
where $\|u\|_1:= \|A(0) u\|$;
\item [$(\mathcal{P}_2)$]
for each $t\geq 0$, the spectrum $\sigma(A(t))$ is contained in the halfplane $\{\zeta \in \C \mid \mathrm{Re}\, \zeta > 0 \}$ and there exists $M>0$ such that, for all $t \geq 0$ and $\zeta \in \C$ with $\mathrm{Re}\, \zeta \leq 0$,
\begin{equation}
\|(\zeta I - A(t))^{-1}\|_{\mathcal{L}(X,X)} \leq
\frac{M}{1+|\zeta|}; \nonumber
\end{equation}
\item[$(\mathcal{P}_3)$]
there exist $L > 0$ and $\gamma \in (0,1]$ such that
\begin{equation}
\|(A(t)-A(s))A(r)^{-1}\|_{\mathcal{L}(X,X)} \leq L
|t-s|^{\gamma}  \ \text{ for all  }\  r,s,t \geq 0. \nonumber
\end{equation}
\end{enumerate}
\begin{remark}{\em
Assumptions $(\mathcal{P}_1)$--$(\mathcal{P}_3)$ are standard and they assure the existence of the parabolic evolution system
determined by the equation $\dot u(t)=-A(t)u(t)$, $t\geq 0$, which we denote by $\{U(t,s)\}_{t \geq s \geq 0}$ (see \cite{Pazy}, \cite{Tanabe} and Proposition \ref{17122016-1936} in Section 2 below). The assumption $(\mathcal{P}_0)$ is natural (and minimal) in the sense that it allows us to consider solutions in the fractional space $X^\alpha$. }
\end{remark}
First we shall consider averaging for linear problems that is
\begin{align}\label{17012017-1620}
\left\{
\begin{array}{ll}
\dot{u}(t)=-A(t/\lma)u(t), \   t>0,\\
u(0)= \bar{u}.
\end{array}
\right.
\end{align}
where $A(t)$, $t\geq 0$, satisfy $(\mathcal{P}_0) - (\mathcal{P}_3)$ and $\bar u \in X^\alpha$. If we define $A^{(\lma)}(t):X^1\to X$, $t\geq 0$, $\lma>0$, by
\begin{equation}\label{20012017-2206}
A^{(\lma)}(t):=A(t/\lma),
\end{equation}
then solutions of \eqref{17012017-1620} are provided by the evolution system
$\{U^{(\lma)}(t,s)\}_{t\geq s\geq 0}$ generated by the operators $A^{(\lma)}(t)$, $t\geq 0$. We shall assume additionally that
\begin{enumerate}
\item[($\mathcal{\widehat{P}}_1)$] there is a constant $K>0$ such that, for any $s,t \geq 0$ with $t>s$ and  $\lma>0$,
\begin{align}
\|U^{(\lma)} (t,s)\|_{\mathcal{L}(X^\alpha, X^\alpha)} \leq K   \ \ \mbox{ and } \ \
\|U^{(\lma)} (t,s)\|_{{\mathcal{L}}(X,X^\alpha)} \leq  K(1+(t-s)^{-\alpha});\nonumber
\end{align}
\end{enumerate}
\begin{enumerate}
\item [$(\mathcal{\widehat{P}}_2)$] there is a linear operator $\widehat A:  X^1 \to X$ generating an analytic $C_0$-semigroup $\{\widehat S(t):X\to X\}_{t\geq 0}$ of bounded linear operators on $X$ such that $D(\widehat A^\alpha)=X^\alpha$ ($\alpha$ is from $({\cal P}_0)$), $\|\widehat S(t)\|_{\mathcal{L}(X,X)} \leq \widehat M$, for all $t \geq 0$ and some constant $\widehat M>0$, and for any $ u \in X^1$
$$
\lim_{\omega \to +\infty} \frac{1}{\omega} \int_0^\omega  A(t) u \d t = \widehat A u.
$$
\end{enumerate}
The main result for the linear problem is as follows.
\begin{theorem}\label{17012017-1853}
Suppose that $(\mathcal{P}_0)$--$(\mathcal{P}_3)$, $(\mathcal{\widehat{P}}_1)$
and $(\mathcal{\widehat{P}}_2)$ hold. Let $(\lma_n)$ in $(0,+\infty)$ be such that $\lma_n \to 0^+$, as $n \to +\infty$. Then
\begin{enumerate}
\item[{\em (i)}] if $(\bar u_n)$ is a sequence in $X^\alpha$ such that $\bar u_n \to \bar u$ in $X$, as $n \to +\infty$, for some $\bar u \in X^\alpha$, then for any $t \geq 0$, $s\in [0,t]$
\begin{equation}
\lim_{n \to +\infty} U^{(\lambda_n)}(t,s) \bar{u}_n= \widehat{S}(t-s) \bar{u} \ \mbox{ in }  \ X^\alpha
\end{equation}
and the convergence  is uniform with respect
to $t\in [0,T]$ and $s\in [0,t-\delta]$, for any fixed $T>0$ and $\delta \in (0,T]$;\\
\item[{\em (ii)}] if $(\bar u_n)$ is a sequence in $X^\alpha$ such that $\bar u_n \to \bar u$ in $X^\alpha$, as $n \to +\infty$, for some $\bar u \in X^\alpha$, then for any $t \geq 0$, $s\in [0,t]$
\begin{equation}
\lim_{n \to +\infty} U^{(\lambda_n)}(t,s) \bar{u}_n= \widehat{S}(t-s) \bar{u} \ \mbox{ in } \ X^\alpha
\end{equation}
and the convergence  is uniform with respect to $t \in [0,T]$  and $s\in [0,t]$ for any fixed $T>0$.
\end{enumerate}
\end{theorem}
\begin{remark}\label{22012017-1350}{\em
(i) The estimates in $(\mathcal{\widehat P}_1)$, where $K$ is independent of $\lma>0$,  do not follow from the construction of parabolic systems for $A^{(\lma)}(t)$, $t\geq 0$, since the constant $L$ in $(\mathcal{P}_3)$ (on which $K$ depends in the Levy construction) in case of $A^{(\lma)}(t)$, $t\geq 0$, is replaced by $L/\lma^\gamma$ (which means that we have no common constant $K$ for the systems $U^{(\lma)}$). However, in applications the estimates $(\mathcal{\widehat P}_1)$ can be derived by different techniques (see Section 6).\\
\indent (ii) The assumption $(\mathcal{\widehat P}_2)$ is satisfied if for instance the coefficients of elliptic operators are periodic or almost periodic in time (see Section 6).}
\end{remark}

Our second main result is for nonlinear evolution problems of the form \eqref{eq:introduction} where the continuous perturbation $F: [0, +\infty)\times X^\alpha \to X$  is supposed to satisfy the following conditions
\begin{enumerate}
\item [$(\mathcal{F}_1)$]  for any $R>0$ there exists  $Q>0$  such that
$$
\|F (t,u_1)-F(t,u_2) \| \leq Q\|u_1-u_2\|_{\alpha} \mbox{ for all } \ t \geq 0 \mbox{ and } u_1, u_2 \in B(0,R);
$$
\item [$(\mathcal{F}_2)$] there is $C>0$ such that $ \| F (t,u) \| \leq C(1+\| u \|_\alpha) \mbox{ for all } t\geq 0 \mbox{ and } u\in X^\alpha$;
\item [$(\widehat{\mathcal{F}})\,$] for each $u \in X^\alpha$ the set $\{F(t, u)\ |\ t \geq
0\}$ is relatively compact and there is a continuous mapping $\widehat F: X^\alpha
\to X$ such that, for any $u \in X^\alpha$,
$$
\widehat F (u) = \lim_{\omega \to +\infty} \frac{1}{\omega} \int_0^\omega F(t,u) \d t.
$$
\end{enumerate}
The following result is a direct consequence of Theorem \ref{06012017-2015}. 
\begin{theorem}\label{20012017-2210} 
Suppose that $(\mathcal{P}_0)$--$(\mathcal{P}_3)$, $(\mathcal{\widehat{P}}_1)$ and $(\mathcal{\widehat{P}}_2)$ as well as  $(\mathcal{F}_1)$,$(\mathcal{F}_2)$ and $(\widehat{\mathcal{F}})$ hold. Let $(\lma_n)$ in $(0, +\infty)$ be such that $\lma_n\to 0^+$, $(\bar u_n)$ be a bounded sequence in $X^\alpha$ and $u_n: [0,
+\infty) \to X^\alpha$ be the solution of $(P_{\lma_n})$ satisfying $u_n(0)= \bar
u_n$, $n \geq 1$. Then
\begin{enumerate}
\item[{\em (i)}] if $\bar u_n \to \bar u_0$ in
$X$ for some $\bar u_0 \in X^\alpha$, then $u_n(t)\to \widehat u(t)$ in $X^\alpha$, where $\widehat u:[0,+\infty) \to X^\alpha$ is the solution of
\begin{equation}  \nonumber
\left\{
 \begin{array}{l}
\dot{u}(t)= \widehat A u(t) + \widehat  F(u(t)), \qquad t >  0,\\
u(0)= \bar u_0,
\end{array}
\right.
\end{equation}
uniformly for $t$ from compact subsets of $(0,+\infty)$;
\item[{\em (ii)}] if $\bar u_n \to \bar u_0$ in
$X^\alpha$ for some $\bar u_0 \in X^\alpha$,
then $u_n(t)\to \widehat u(t)$ in $X^\alpha$ uniformly for $t$ from compact subsets of  $[0,+\infty)$.
\end{enumerate}
\end{theorem}
\begin{remark}{\em
(i) Conditions $({\cal F}_1)$ and $({\cal F}_2)$ assure the existence of the unique  (mild) solution $u\in C([0,+\infty), X^\alpha)$ for $(P_\lma)$ with the initial condition $u(0)=\bar u$, for any $\bar u\in X^\alpha$ (see Section 4).\\
\indent (ii) Condition $({\cal\widehat{F}})$ is satisfied if for instance $F$ is $T$-periodic in time with $\widehat F:X^\alpha\to X$ given by $\widehat F(u)=\frac{1}{T}\int_{0}^{T}F(t,u)\d t$, $u\in X^\alpha$.
}\end{remark}
\indent The paper is organized as follows. In Section 2 we briefly recall the
classical result concerning the existence of evolution systems in the parabolic case
and provide the estimates comparing two evolutions systems determined by two different time dependent families of sectorial operators.
Section 3 is devoted to averaging technique for linear abstract parabolic
equations. In section 4 we deal with the continuity properties with respect to $A(t)$'s and $F$. Section 5 shows that averaging principle holds for nonlinear parabolic evolution problems. Finally Section 6 provides an example of application of our results to parabolic partial differential equations with time dependent coefficients.

\section{Parabolic evolution systems}

A family of bounded linear operators $\{U(t,s): X \to X\}_{0 \leq s
\leq t \leq T}$, $T>0$, on a Banach space $X$ is called an
\emph{evolution system} provided $U(t,t)=I$ for each $t \in [0,T]$,
$U(t,s)=U(t,r)U(r,s)$ whenever $0 \leq s \leq r \leq t \leq T$ and for any
$\bar{u} \in X$, the map $(s,t) \mapsto U(t,s)\bar{u}$ is
continuous. Evolution systems appear as translation along trajectories operators of differential equations that involve
families of operators which are dependent on time. \\
\indent Consider the problem
\begin{align}\label{17122016-1914}
\left\{
\begin{array}{ll}
\dot{u}(t)=-A(t)u(t), \ \  t \in (s,T]\\
u(s)= \bar{u},
\end{array}
\right.
\end{align}
where $T>0$, $\bar u \in X$ is arbitrary and the family $\{A(t)\}_{t \in [0,T]}$ of linear operators in $X$ satisfies the parabolic conditions $(\mathcal{P}_1)-(\mathcal{P}_3)$ (restricted to $s,t\in [0,T]$).
\begin{remark} \label{23022016-1024}
\em{ The assumptions $(\mathcal{P}_1)$ and $(\mathcal{P}_2)$ state that for each $t
\in [0,T]$, the operator $A(t): D(A(t)) \to X$ is sectorial, which is equivalent to the fact that $-A(t)$ generates an analytic  $C_0$-semigroup -- see \cite{Pazy}.}
\end{remark}
Recall that a function $u\in C([s,T],X)\cap C^1((s,T])$ is called a {\em classical solution} of the problem (\ref{17122016-1914})
if $u(t) \in D(A(t))$ for $ t\in (s,T]$, $u$ is continuously
differentiable on $(s,T]$ and satisfies (\ref{17122016-1914}).
It is known that, under assumptions $(\mathcal{P}_1)-(\mathcal{P}_3)$,
there exists a unique classical solution $u=u_{s, \bar u}: [s ,T) \to X$ of
$(\ref{17122016-1914})$. The corresponding parabolic evolution system determined by the family of operators $\{A(t)\}_{t \geq 0}$ can be constructed by use of the the Levi method and then $u_{s,\bar u} (t)= U(t,s)\bar{u}$. For details we refer to \cite{Pazy} or \cite{Tanabe} where the following result comes from.
\begin{proposition}[see \cite{Pazy}, Ch. 5, Th. 6.1] \label{17122016-1936}
Let $\{A(t)\}_{t \in [0,T]}$ be a family of linear operators in a Banach space $X$ satisfying conditions
$(\mathcal{P}_1)-(\mathcal{P}_3)$. Then there exists a unique
evolution system $\{U(t,s):X \to X\}_{0 \leq s \leq t \leq T}$
with the following properties
\begin{enumerate}
\item[\emph{(i)}] there exists a constant $D_1 > 0$ such that $\|U(t,s)\|_{\mathcal{L}(X,X)} \leq D_1$ whenever
$0 \leq s \leq t\leq T$;
\item[\emph{(ii)}] $U(t,s)X \subset X^1$ if $0 \leq s < t \leq T$, for each $u\in X$, the function $t\mapsto U(t,s)u$
is differentiable in $X$,  $\frac{\partial}{\partial t}U(t,s)\in \mathcal{L}(X,X)$,
the mapping $(s,t)\mapsto \frac{\partial}{\partial t}U(t,s)$ is strongly continuous on the set $\{ (s,t)\in [0,T]\times [0,T] \mid s < t \}$ and 
\vspace{-4pt}
\begin{align}
& & \frac{\partial}{\partial t}U(t,s) + A(t)U(t,s)\  & =  \  0 &\quad \text{if\ } 0 \leq s < t \leq T, \nonumber \\
& & \left\|\frac{\partial}{\partial t} U(t,s) \right\|_{\mathcal{L}(X,X)} &\leq \frac{D_2}{t-s} & \quad  \text{ if } 0 \leq s < t \leq T, \nonumber \\
& & \| A(t)U(t,s)A(s)^{-1} \|_{\mathcal{L}(X,X)} \ &\leq  \ D_3 & \quad
\text{ if  } 0 \leq s \leq t \leq T \nonumber
\end{align}
for some constants $D_2, D_3>0$;
\item[\emph{(iii)}] for every $\bar{v} \in X^1$ and $t \in (0,T]$, $[0,t]\ni s\mapsto U(t,s)\bar{v} \in X$ is
differentiable and
\begin{equation}
\frac{\partial}{\partial s}U(t,s)\bar{v} =
U(t,s)A(s)\bar{v}.\nonumber
\end{equation}
\end{enumerate}
\end{proposition}
\begin{remark}\label{11102016-0709}
{\em 
 (i)  If we consider $\{ A(t)\}_{t\geq 0}$ satisfying
$(\mathcal{P}_1)-(\mathcal{P}_3)$ (for all $t\geq 0$), then we obtain
	$\{ U(t,s) \}_{t\geq s\geq 0}$ by increasing $T$. Nevertheless, the constants $D_1$, $D_2$ and  $D_3$ may depend on $T$. Therefore we perform estimations with use of the inequalities involving those constants only for $s$ and $t$ from bounded intervals.\\
\indent (ii) It is clear that $({\cal P}_1)$ implies $\|A(t)\|_{{\cal L}(X^1,X)} \leq C_2$ for all $t\geq 0$.\\
\indent (iii) It comes from $({\cal P}_3)$ that
$[0,+\infty)\ni t\mapsto A(t)\in {\cal L}(X^1,X)$ is continuous. Indeed, for any $v\in X^1$ and $t,s\geq 0$, one has
$$
\|(A(t)-A(s))v\| \leq \|(A(t)-A(s)) A(0)^{-1} A(0)v\|
\leq L|t-s|^\gamma \|v\|_1.
$$}
\end{remark}
\begin{corollary}\label{12042017-1251}
For any sequences
$(t_n)$, $(s_n)$ in $[0,+\infty)$ and $(\bar u_n)$ in $X$ such that
$t_n\to t$, $s_n\to s$ and $\bar u_n\to \bar u$ for some $t,s$ with $t>s$ 
and $\bar u\in X$
$$
U(t_n,s_n)\bar u_n \to U(t,s)\bar u \ \mbox{ in } \ X^1, \mbox{ as } n\to +\infty.
$$
\end{corollary}
\noindent {\bf Proof:} Observe that for large $n$
$$
\|U(t_n,s_n)\bar u_n - U(t,s)\bar u\|_1 \leq C_1^{-1} (\alpha_n +\beta_n + \gamma_n)
$$
where
\begin{eqnarray*}
 \alpha_n & := & \|A(t_n)U(t_n, s_n) (\bar u_n - \bar u)\|,\\
 \beta_n & := & \|A(t_n)U(t_n,s_n)\bar u- A(t)U(t,s)\bar u \|,\\
 \gamma_n& := & \| (A(t)-A(t_n) U(t,s)\bar u\|. 
\end{eqnarray*}
Clearly, due to Proposition \ref{17122016-1936} (ii), for large $n$ one has $\alpha_n \leq D_2 \|\bar u_n-\bar u \|/(t_n-s_n)$, which means that  $\alpha_n \to 0$. Since $(\tau,\rho)\mapsto\frac{\part U}{\part \tau}(\tau,\rho)\bar u\in X$ is continuous, we see that $\beta_n \to 0$. Finally, in view of Remark \ref{11102016-0709} (iii),
$\gamma_n \leq L|t-t_n|^\gamma \|U(t,s)\bar u \|_1 \to 0$. \hfill $\square$\\

\begin{remark} {\em 
Hence, if we assume additionally $(\mathcal{P}_0)$, then it makes sense to consider $U(t,s):X^\alpha \to X^\alpha$ for $s,t\geq 0$ with $s\leq t$ and we see that $U(t_n, s_n) \bar u_n \to U(t,s)\bar u$ in $X^\alpha$, whenever $(t_n)$, $(s_n)$, $(\bar u_n)$ are as in Corollary \ref{12042017-1251}.
}
\end{remark}

\section{Continuity properties and proof of Theorem \ref{17012017-1853}}

Consider families of linear operators $\{A_n(t)\}_{t \geq 0}$, $n \geq 1$, satisfying $({\cal P}_0)-({\cal P}_3)$ with common spaces $X^1$ and $X^\alpha$, $\alpha\in (0,1)$, constants $C_{1,\alpha}, C_{2,\alpha}, C_1, C_2$ and $M$ (not necessarily $L$). Denote by $\{U_n (t,s)\}_{t \geq s \geq 0}$ the evolution system determined by $\{A_n(t)\}_{t\geq 0}$, $n\geq 1$, and assume additionally that the following conditions hold:
\begin{enumerate}
	\item[$(\widetilde{\cal P}_1)$] there exists $K>0$ such that, for any $s,t \geq 0$ with $t>s$ and $n\geq 1$,
	\begin{equation}
	\|U_n(t,s)\|_{{\cal L}(X^\alpha,X^\alpha)} \leq K \ \ \mbox{ and } \ \ 
	\|U_n(t,s)\|_{{\cal L}(X,X^\alpha)}\leq K(1+(t-s)^{-\alpha}); \nonumber
	\end{equation}
\item[$(\widetilde{\cal P}_2)$] there $\{S_0(t)\}_{t\geq 0}$ is an analytic $C_0$-semigroup generated by $-A_0$, where $A_0:X^1 \to X$, such that for any $\beta \geq 0$ and $t>0$, $\|A_0^\beta S_0(t)\|_{\mathcal{L}(X,X)}\leq M_\beta t^{-\beta}$ for some $M_\beta>0$ and, for each $t \geq 0$,
\begin{equation}
\lim_{n \to +\infty} \int_0^t \left( A_n(r)-A_0 \right) \d r = 0 \ \textnormal{ in }\ {\cal L}(X^1, X) \nonumber
\end{equation}
and the convergence is uniform with respect to $t \in [0,T]$,  for any fixed $T>0$.
\end{enumerate}
Our aim is the following continuity property.
\begin{theorem}\label{20042017-1423}
Under the above assumptions suppose additionally that
\begin{enumerate}
\item[{\em (i)}] if $(\bar u_n)$ is a sequence in $X^\alpha$ such that $\bar u_n \to \bar u$ in $X$, as $n \to +\infty$, for some $\bar u \in X^\alpha$, then for any $t \geq 0$, $s\in [0,t]$
\begin{equation}
\lim_{n \to +\infty} U_n (t,s) \bar{u}_n= S_0(t-s) \bar{u} \ \mbox{ in }  \ X^\alpha
\end{equation}
and the convergence  is uniform with respect
to $t\in [\delta,T]$ and $s\in [0,t-\delta]$, for any fixed $T>0$ and $\delta \in (0,T)$;\\
\item[{\em (ii)}] if $(\bar u_n)$ is a sequence in $X^\alpha$ such that $\bar u_n \to \bar u$ in $X^\alpha$, as $n \to +\infty$, for some $\bar u \in X^\alpha$, then for any $t \geq 0$, $s\in [0,t]$
\begin{equation}
\lim_{n \to +\infty} U_{n}(t,s) \bar{u}_n= S_0 (t-s) \bar{u} \ \mbox{ in } \ X^\alpha
\end{equation}
and the convergence  is uniform with respect to $t \in [0,T]$  and $s\in [0,t]$ for any fixed $T>0$.
\end{enumerate}
\end{theorem}
\noindent We shall need some preparation before proving the theorem.
\begin{lemma}\label{11032017-1828}
Under the above assumptions, for any $T>0$, $t \in [0,T]$, $s \in [0,t]$, $\bar u \in X^1$ and $n \geq 1$, 
\end{lemma}
\begin{equation}
S_0(t-s)\bar u - U_n(t,s)\bar u = \int_s^t U_n(t,r) \left( A_n(r)-A_0\right) S_0(r-s) \bar u \d r. \nonumber
\end{equation}
\noindent {\bf Proof:} For given $\bar u \in X^1$, $s,t \in [0,T]$ with $s<t$ define $\phi: [s,t] \to X$ by
$$
\phi(r):= U_n(t,r)S_0(r-s) \bar u.
$$
Due to Proposition \ref{17122016-1936}, the map $\phi$ is differentiable in $X$ on the interval $(s,t)$ and we have
\begin{align}
\phi'(r) &= U_n(t,r)A_n(r)S_0(r-s) \bar u - U_n(t,r)A_0 S_0(r-s)\bar u \nonumber \\
& = U_n(t,r) \left( A_n(r)-A_0\right) S_0(r-s) \bar u, \label{09032017-2147}
\end{align}
whenever $0<s<r<t \leq T$. Therefore one has
\begin{equation}
S_0(t-s) \bar u - U_n(t,s)\bar u = \phi(t)-\phi(s) = \int_s^t \phi'(r) \d r, \nonumber
\end{equation}
which together with \eqref{09032017-2147} 
yields the assertion. \hfill $\square$\\

\begin{lemma}\label{12032017-2243}
Suppose that $F_n:[0,+\infty) \to \mathcal{L}(X,X)$, $n \geq 1$, are continuous functions such that for all $t \geq 0$ and $n\geq 1$, $\|F_n(t)\|_{\mathcal{L}(X,X)} \leq R_F$  for some $R_F>0$ and
\begin{equation}
\lim_{n \to +\infty} \int_0^t F_n(r) \d r = 0 \textnormal{ in } X, \label{12032017-1931}
\end{equation} 
uniformly with respect to $t \in [0,T]$ for any $T>0$. Then 
\begin{equation}
\lim_{n \to +\infty} \int_s^t F_n(r) S_0(r-s) \bar u \d r = 0 \textnormal{ in }X,\nonumber
\end{equation}
uniformly with respect to $\bar u$ from bounded sets in $X$, $t \in [0,T]$ and $s \in [0,t]$.
\end{lemma}
\noindent {\bf Proof:} Fix $T>0$ and $R>0$. Take an arbitrary $\varepsilon>0$ and $\eta > (0,T)$ such that $\eta M_0 R R_F<\varepsilon/3$. Clearly, if $t\in [0,T]$ and $0\leq t-s<\eta$, then 
\begin{equation}\label{08052017-1943}
\left\| \int_s^t F_n(r) S_0(r-s) \bar u \d r \right\| \leq
(t-s)M_0 R R_F < \varepsilon/3 <\varepsilon
\end{equation}
for all $n\geq 1$ and $\bar u\in X$ with $\|\bar u\|\leq R$. And, for $t \in [\eta,T]$, $s\geq 0$ such that $t-s \geq \eta$ and $\bar u \in X$ with $\|\bar u\|\leq R$, one gets
\begin{equation}\label{25042017-1532}
\int_s^t F_n(r)S_0(r-s) \bar u \d r = I_{1,n}(\bar u, s, t)+I_{2,n}(\bar u, s, t),
\end{equation} 
where
\begin{align}
I_{1,n} (\bar u, s, t)&:= \int_s^{s+\eta} F_n(r)S_0(r-s) \bar u \d r, \nonumber \\
I_{2,n}(\bar u, s, t)&:= \int_{s+\eta}^t F_n(r)S_0(r-s) \bar u \d r. \nonumber
\end{align}
As before we note that, for all $n \geq 1$,
\begin{equation}\label{12032017-1934}
\|I_{1,n}(\bar u, s,t)\| \leq  \eta M_0 R R_F < \varepsilon/3. 
\end{equation}
Further, integrating by parts we get
\begin{align}
I_{2,n} (\bar u,s,t)& = \int_{s+\eta}^t \frac{\d}{\d r}\left(\int_{s+\eta}^r F_n(\tau)\d \tau \right) S_0(r-s)\bar u \d r \nonumber \\
& = \left[\int_{s+\eta}^t F_n(\tau)\d \tau \right] S_0(t-s)\bar u + \int_{s+\eta}^t \left(\left[\int_{s+\eta}^r F_n(\tau)\d \tau\right]  A_0 S_0(r-s)\bar u \right) \d r.\nonumber
\end{align}
By \eqref{12032017-1931} we find $N \geq 1$ such that for any $n \geq N$ and
$t\in [0,T]$,
$$
\left\|\int_0^t F_n(r) \d r\right \| \leq \frac{\varepsilon}{6R}\cdot\min\{1/M_0, \eta/M_1 T\}.
$$
Hence, for all $n \geq N$, $t\in [\eta,T]$ and $s\in [0,t-\eta]$,  
\begin{align}
\|I_{2,n}(\bar u,s,t)\|&\leq M_0 R \left\|\int_{s+\eta}^t F_n(\tau)\d \tau \right\| + \int_{s+\eta}^t \frac{M_1 R}{r-s} \left\|\int_{s+\eta}^r F_n(\tau)\d \tau \right\| \d r \nonumber \\
& \leq M_0R \left\|\int_{s+\eta}^{t} F_n(\tau)\d \tau\right\| + \frac{M_1 R T}{\eta}\left\|\int_{s+\eta}^{t} F_n(\tau)\d \tau\right\|\leq \varepsilon/3+\varepsilon/3=2\varepsilon/3,\nonumber 
\end{align}
which together with \eqref{25042017-1532} and \eqref{12032017-1934}  implies that, for any $n \geq N$,
$$
\left\|\int_s^t F_n(r)S_0(r-s) \bar u \d r \right\|\leq \varepsilon 
\ \mbox{ for all } \ t\in [\eta,T] \mbox{ and } s\in [0,t-\eta].
$$
This along with \eqref{08052017-1943} gives the assertion. \hfill $\square$\\

\begin{lemma} \label{12032017-2258}
Let $T>0$ and assume that functions $G_n: \Delta  \to \mathcal{L}(X,X)$, $n\geq 1$, where $\Delta:=\{(t,s) \in [0,T]\times [0,T]\ |\ t\geq s \}$ are such that $\|G_n(t,s)\|_{\mathcal{L}(X,X)}\leq R_G$ for all $(t,s) \in \Delta$, $n \geq 1$ and some $R_G>0$, the mapping $(t,s) \mapsto G_n(t,s) \bar u$ is continuous for each $\bar u \in X$ and $n \geq 1$
\begin{equation}
\lim_{n \to +\infty} \int_s^t G_n(r,s) \bar u \d r = 0 \ \textnormal{ in } \ X,\label{12032017-1649}
\end{equation}
uniformly with respect to $\bar u$ from bounded sets in $X$, $t \in [0,T]$ and $s\in [0, t]$. 
Then
\begin{equation}
\lim_{n \to +\infty} \int_s^t S_0(t-r) G_n(r,s) \bar u \d r =0 \ \textnormal{ in } \ X^\alpha, \nonumber
\end{equation}
and the convergence is uniform with respect to $\bar u$ from bounded sets in $X$, $t \in [0,T]$ and $s \in [0,t]$.
\end{lemma}
\noindent {\bf Proof:} Fix an arbitrary $R>0$ and take any $\varepsilon>0$.  Let $\eta \in (0,T)$ be such that
$$
M_\alpha R_G R  \frac{\eta^{1-\alpha}}{1-\alpha}<\varepsilon/3.
$$ 
Clearly, if $t\in [0,T]$ and $0<t-s<\eta$, then
\begin{equation}\label{08052017-1956}
\left\| \int_s^t S_0(t-r)G_n(r,s) \bar u \d r \right\|_\alpha
\leq M_\alpha R_G R \int_{s}^{t} (t-s)^{-\alpha} \d r  \leq M_\alpha R_G R  \frac{\eta^{1-\alpha}}{1-\alpha}< \varepsilon
\end{equation}
for all $n\geq 1$ and $\bar u\in X$ with $\|\bar u\|\leq R$.
And, if $t \in [\eta,T]$, $t-s\geq\eta$ and $\bar u \in X$ with $\|\bar u  \|\leq R$, then
\begin{equation}\label{27042017-1837}
\int_s^t S_0(t-r)G_n(r,s) \bar u \d r = J_{1,n}(\bar u, s,t) + J_{2,n} (\bar u, s,t)
\end{equation}
where
\begin{align*}
J_{1,n} (\bar u, s,t) &:= \int_{t-\eta}^t S_0(t-r)G_n(r,s) \bar u \d r, \\\
J_{2,n}(\bar u, s,t) &:= \int_s^{t-\eta} S_0(t-r)G_n(r,s) \bar u \d r.  
\end{align*}
For all $n \geq 1$, we get
\begin{eqnarray}\nonumber
\|J_{1,n} (\bar u, s,t)\|_\alpha & \leq & \int_{t-\eta}^t \|S_0(t-r)\|_{\mathcal{L}(X,X^\alpha)}\|G_n(r,s)\|_{\mathcal{L}(X,X)}\|\bar u\|\d r \\
& \leq & M_\alpha R_G R  \frac{\eta^{1-\alpha}}{1-\alpha}<\varepsilon/3. \label{12032017-1656}
\end{eqnarray}
Now observe that for all $n \geq 1$ one has
\begin{align}
J_{2,n} (\bar u, s,t)&= \int_s^{t-\eta} S_0(t-r) \frac{\d}{\d r}\!\left(-\int_r^{t-\eta} \!G_n(\tau,s)\bar u \d \tau \right) \d r \nonumber \\
& = S_0(t-s) \left[\int_s^{t-\eta} \!\!\!\! G_n(\tau, s) \bar u \d \tau\right] +\!\int_s^{t-\eta} \!\!\!\! A_0 S_0(t-r) \left[\int_r^{t-\eta}\!\!\!\! G_n(\tau,s) \bar u \d \tau\right]\!\!\d r.\nonumber
\end{align} 
Obviously, by \eqref{12032017-1649} it follows that there is $N\geq 1$ such that for all $n\geq N$, $t\in [\delta, T]$, $s\in [0,t-\delta]$ and $\bar u\in X$ with $\|\bar u\|\leq R$, we obtain 
$$
\left\|\int_s^t G_n(r,s) \bar u \d r\right\| \leq \frac{\varepsilon}{3}\cdot\min\{\eta^\alpha/ M_\alpha,  \eta^{1+\alpha}/ M_{1+\alpha} T\}.
$$
This gives, for all $n\geq N$, $t\in [\eta, T]$, $s\in [0,t-\eta]$ and $\bar u\in X$ with $\|\bar u\|\leq R$, 
\begin{align*}
\|J_{2,n}(\bar u,s,t)\|_\alpha &\leq \frac{M_\alpha}{(t-s)^\alpha} \left\|\int_s^{t-\eta}\!\!\!\! G_n(\tau,s)\bar u \d \tau \right\| +\int_s^{t-\eta} \frac{M_{1+\alpha}}{(t-r)^{1+\alpha}} \left\|\int_r^{t-\eta}\!\!\!\! G_n(\tau,s)\bar u \d \tau \right\| \d r \nonumber \\
&\leq \frac{M_\alpha}{\eta^\alpha} \left\|\int_s^{t-\eta} \!\!\!\! G_n(\tau,s)\bar u \d \tau \right\|+  \frac{M_{1+\alpha}T}{\eta^{1+\alpha}} \cdot \sup_{r\in [s,t-\eta]}\left\|\int_r^{t-\eta} \!\!\!\! G_n(\tau,s)\bar u \d \tau \right\|  \nonumber \\
&\leq  \varepsilon/3+\varepsilon/3=2\varepsilon/3,\nonumber
\end{align*}
which together with \eqref{27042017-1837}, \eqref{12032017-1656} as well as \eqref{08052017-1956} yields the assertion.  
\hfill $\square$\\
\begin{lemma}\label{11032017-1914}
For any $T>0$ 
\begin{equation*}
\int_s^t S_0(t-r) \left( A_n(r)-A_0 \right) S_0(r-s) \bar u \d r  \to 0 \ \textnormal{ in } \ X^\alpha, 
\end{equation*}
uniformly with respect to $\bar u$ from bounded sets in $X^1$, $t \in [0,T]$ and $s \in [0, t]$.
\end{lemma}
\noindent {\bf Proof:} Fix $T>0$. Define $F_n:[0,+\infty) \to \mathcal{L}(X,X)$, $n \geq 1$, by 
$$
F_n(t):= (A_n(t)-A_0)A_0^{-1}
$$
and let $G_n:\Delta \to \mathcal{L}(X,X)$, $n\geq 1$, be given by
$$
G_n(t,s):= F_n(t)S_0(t-s).
$$
Clearly, for any $t \in [0,T]$, $s \in [0,t]$, $n \geq 1$ and $\bar u \in X^1$, one has
\begin{equation}\label{12032017-2233}
 S_0(t-r)(A_n(r)-A_0)S_0(r-s)\bar u =  S_0(t-r)G_n(r,s)A_0 \bar u.
\end{equation}
Observe that, for each $n \geq 1$, $F_n$ is continuous and  it follows from $({\cal P}_1)$ that $\|F_n(t)\|_{\mathcal{L}(X,X)} \leq 2C_2/C_1$, for $t \geq 0$ and $n\geq 1$. Next, in view of $(\widetilde{\cal P}_2)$, we infer that, for any $\bar v \in X$,
$$
\left\| \int_0^t F_n(r) \bar v \d r\right\| \leq  \left\| 
\int_0^t (A_n(r)-A_0) \d r\right\|\cdot \|A_0^{-1}\|_{{\cal L}(X,X^1)} \to 0\ \textnormal{ in } \ X,
$$ 
uniformly with respect to $t \in [0,T]$. Thus, using  Lemma \ref{12032017-2243} we conclude that, for $\bar v\in X$,
\begin{equation}\nonumber
\lim_{n \to +\infty} \int_s^t G_n(r,s) \bar v \d r = 0 \textnormal{ in } X
\end{equation}
and the convergence is uniform with respect to $\bar v$ from bounded sets in $X$, $t \in [0,T]$ and $s \in [0, t]$. Further, note that for each $n \geq 1$ and $\bar v \in X$, the mapping $\Delta \ni (t,s) \mapsto G_n(t,s)\bar v$ is continuous and, in view of $(\widetilde{\cal P}_2)$, for all $n \geq 1$ and $(t,s) \in \Delta$,
$$
\|G_n(t,s)\|_{\mathcal{L}(X,X)}\leq \|F_n(t)\|_{\mathcal{L}(X,X)}\|S_0(t-s)\|_{\mathcal{L}(X,X)}
\leq 2C_2 M_0/C_1.
$$
Hence, applying Lemma \ref{12032017-2258}, we get
\begin{equation}\nonumber
\lim_{n \to +\infty} \int_s^t S_0(t-r)G_n(r,s) \bar v \d r =0 \textnormal{ in } X^\alpha,
\end{equation}
uniformly with respect to $\bar v$ from bounded sets in $X$, $t\in [0,T]$ and $s\in [0,t]$ . This together with \eqref{12032017-2233} ends the proof, since  $\|A_0 \bar u\| \leq C_2 \|u\|_1$ for any $\bar u \in X^1$.
\hfill $\square$\\

\noindent \textbf{Proof of Theorem \ref{20042017-1423}:} \\
Fix $T>0$. The assertion (ii) will follows if we prove that, for each $t\in [0,T]$  and $s\in [0,t]$,
\begin{equation}\label{11032017-1806}
\lim_{n \to +\infty} \| U_n(t,s) - S_0(t-s) \|_{{\cal L}(X^1,X^\alpha)} =0
\end{equation}
and the convergence is uniform with respect to $t \in [0,T]$ and $s \in [0,t]$. 
To this end we use Lemma \ref{11032017-1828} to see that, for any $\bar u \in X^1$ with $\|\bar u\|_1=1$, $t\in [0,T]$, $s \in [0,t]$ and $n \geq 1$, 
\begin{align*}
S_0(t-s)\bar u - U_n(t,s)\bar u & =  \int_s^t U_n(t,r)\left(A_n(r)-A_0\right) S_0(r-s)\bar u \d r \\
& =   K_{1,n}(s,t)+K_{2,n}(s,t)
\end{align*} 
where
\begin{align}
K_{1,n}(s,t)&:= \int_s^t S_0 (t-r)\left(A_n(r)-A_0 \right)S_0 (r-s)\bar u \d r, \nonumber \\
K_{2,n}(s,t)&:= \int_s^t \left(U_n(t-r)-S_0(t-r)\right)\left(A_n(r)-A_0\right)S_0(r-s)\bar u \d r. \nonumber
\end{align}
Clearly, for $t\in [0,T]$ and $s\in [0,t]$, $\|K_{1,n}(s,t)\|_{\alpha}\leq a_n(T)$ with 
\begin{equation}
a_n (T):= \sup_{\bar w \in X^1,\ \|\bar w \|_1 \leq 1,\  t\in [0,T],\ s \in [0,t]}\left\| \int_s^t S_0(t-r)\left(A_n(r)-A_0 \right)S_0(r-s)\bar w \d r \right\|_\alpha. \nonumber
\end{equation}
Further, in view of Remark \ref{11102016-0709} (ii) and $(\widetilde{\cal P}_1)$, we get
\begin{align}
\|K_{2,n}(s,t)\|_\alpha &\leq \int_s^t \| U_n(t,r)- S_0(t-r)\|_{\mathcal{L}(X,X^\alpha)} \|A_n(r)- A_0\|_{\mathcal{L}(X^1,X)} \|S_0 (t-r)\bar u\|_1 \d r \nonumber \\
& \leq 2 \tilde C C_2 M_0  \int_{s}^{t} \|  S_0(t-r)- U_n(t,r)\|_{\mathcal{L}(X^1,X^\alpha)} \d r \nonumber
\end{align} 
where $\tilde C >0$ is the constant related to the embedding $X^1 \hookrightarrow X$. Hence, summing up, for any $t \in [0,T]$ and $s \in [0,t]$,
\begin{equation*}
\|S_0 (t-s) - U_n(t,s)\|_{\mathcal{L}(X^1,X^\alpha)} \leq a_n (T) + 2\tilde C C_2 M_0 \int_{s}^{t} \|  S_0 (t-r)-U_n(t,r)\|_{\mathcal{L}(X^1,X^\alpha)} \d r. 
\end{equation*}
By the Gronwall inequality, we obtain for all $t\in [0, T]$ and $s\in [0,t]$,
\begin{equation}\label{08052017-1721}
\|S_0 (t-s) - U_n(t,s)\|_{\mathcal{L}(X^1,X^\alpha)} \leq  a_n(T) \exp \left(2\tilde C C_2 M_0 T \right).
\end{equation}
According to Lemma \ref{11032017-1914}, $a_n(T) \to 0$ as $n\to +\infty$, which along with \eqref{08052017-1721} proves \eqref{11032017-1806}.\\ 
\indent In order to show (i), fix $\delta\in (0,T)$ and let $(\bar u_n)$ be a sequence in $X^\alpha$ such that $\bar u_n \to \bar u$ in $X$, as $n \to +\infty$ for some $\bar u \in X^\alpha$. For any $\varepsilon>0$ we need to find $N \geq 1$ such that, for all $n \geq N$, $t\in [\delta,T]$ and $s\in [0,t-\delta]$,
$$
\| U_n (t,s)\bar u_n - S_0 (t-s)\bar u\|_\alpha <\varepsilon.
$$ 
To this end note that, by the density of $X^1$ in $X^\alpha$, we can choose $\bar u_0\in X^1$ with $\|\bar u_0\|_1=1$ such that
$$
\|\bar u-\bar u_0\|_\alpha <\min \left\{\eta/C_\alpha, \varepsilon/(3  M_0)\right\}
$$
where $\eta:= \varepsilon/\left(6K(1+\delta^{-\alpha})\right)$ and $C_\alpha>0$ is the constant in the inequality $\|u\| \leq C_\alpha \|u\|_\alpha$ for $u \in X^\alpha$. In view of \eqref{11032017-1806}, we get $N_0 \geq 1$ such that
\begin{equation}\label{07022017-0406}
\| U_n(t,s)\bar u_0 -  S_0(t-s)\bar u_0 \|_\alpha <\varepsilon/3 \ \ \text{ for all } \ \ n \geq N_0,\, t\in [0,T],\, s\in [0,t],
\end{equation}
Let $N_1 \geq 1$ be such that, for all $n \geq N_1$, $\|\bar u_n - \bar u\| < \eta$.
This, together with $(\mathcal{\widetilde {P}}_1)$, \eqref{07022017-0406} and $(\mathcal{\widetilde P}_2)$ \footnote{In view of $(\mathcal{\widetilde P}_2)$, $\| S_0(t)\|_{\mathcal{L}(X^\alpha,X^\alpha)} \leq M_0$ for all $t \geq 0$.}, implies that for $n \geq N$, where $N:=\max\{N_0, N_1\}$, $t\in [\delta,T]$ and $s\in [0,t-\delta]$, we have
\begin{align*}
\| U_n(t,s)\bar u_n - S_0 (t-s)\bar u\|_\alpha  & \leq     
\| U_n(t,s)(\bar u_n -\bar u_0)\|_\alpha + \|U_n(t,s)\bar u_0- S_0 (t-s)\bar u_0 \|_\alpha\\
 &   +  \|S_0 (t-s) (\bar u_0 -\bar u) \|_\alpha\\
& <     K(1+(t-s)^{-\alpha})\|\bar u_n - \bar u_0\|  + \varepsilon /3 + \\
&  + \| S_0 (t-s) \|_{{\cal L}(X^\alpha, X^\alpha)} \|\bar u_0 - \bar u\|_\alpha\\
& < K(1+\delta^{-\alpha})(\|\bar u_n - \bar u\|+C_\alpha\|\bar u - \bar u_0 \|_\alpha) + \varepsilon/3 + M_0 (\epsilon /3 M_0) \\
&  < K(1+\delta^{-\alpha})(\eta+C_\alpha (\eta/C_\alpha)) + \epsilon/3 + \epsilon /3 = \epsilon,
\end{align*}
which ends the proof of (ii). \hfill $\square$\\

\noindent \textbf{Proof of Theorem \ref{17012017-1853}:} 
Take any $(\lma_n)$ such that $\lma_n \to 0^+$ as $n\to +\infty$. 
Put $A_n(t):=A^{(\lma_n)}(t)$ and $U_n(t,s):=U^{(\lma_n)}(t,s)$ and observe that the families of operators $\{A_n(t)\}_{t \geq 0}$, $n \geq 1$,
satisfy conditions $(\mathcal{P}_0)$--$(\mathcal{P}_3)$ with the same constants $C_{1,\alpha}, C_{2,\alpha},C_1, C_2, M>0$ (and the constants $L/\lma_n^\gamma$ instead of $L$). It is immediate that $({\cal \widehat P}_1)$ implies $({\cal \widetilde P}_1)$. Let $A_0 := \widehat A$ and $S_0 (t):=\widehat S (t)$. To show that $({\cal \widehat P}_2)$ implies
$({\cal \widetilde P}_2)$ note that
\begin{align*}
\int_0^t (A_n(r)-\widehat A) \d r& =\int_0^t (A(r/\lma_n)-\widehat A) \d r = \lma_n \int_0^{t/\lma_n} (A(\rho)-\widehat{A})  \d \rho \\
& = t \cdot \frac{1}{t/\lma_n}\int_0^{t/\lma_n} ( A(\rho) - \widehat{A})  \d \rho.
\end{align*}
Take $\delta \in (0,T)$ and observe that for $t \in [0,\delta]$, by Remark \ref{11102016-0709} (ii) we get
\begin{equation}\label{15032017-1343}
\left\|\int_{0}^{t} (A_n (r)-A_0) \d r\right\|_{{\cal L}(X^1,X)} \leq 2 C_2 \delta,
\end{equation}
and, for $t \in [\delta,T]$, one has 
\begin{equation}
\left\| \int_0^t (A_n(r)- A_0) \d r \right\|_{{\cal L}(X^1,X)}  \leq T \cdot \left\| \frac{1}{t/\lma_n}  \int_0^{t/\lma_n} ( A(\rho) - \widehat{A})\d\rho\right\|_{{\cal L}(X^1,X)},   \nonumber
\end{equation}
which in view of $(\mathcal{\widehat{P}}_2)$ and together with \eqref{15032017-1343} yields $({\cal \widetilde P}_2)$. We obtain the assertion by applying Theorem \ref{20042017-1423}. \hfill $\square$\\

\section{Continuity of nonlinear evolution systems}

We shall consider the following abstract parabolic equation
\begin{equation}  \label{07022014-1007}
\left\{   \begin{array}{l}
\dot{u}(t)= -A(t) u(t) + F(t, u(t)), \qquad t >  0,\\
u(0)= \bar u,
\end{array}  \right.
\end{equation}
where $\bar u\in X^\alpha$, the family $\{A(t): D(A(t)) \to X\}_{t \geq 0}$ of linear operators on a Banach space $X$ satisfies $({\cal P}_0)$--$({\cal P}_3)$ and $F:[0, +\infty) \times X^\alpha \to X$ is a continuous mapping such that $(\mathcal{F}_1)$ and $(\mathcal{F}2)$ hold.\\
\indent We shall say that a continuous $u:[0,T) \to X^\alpha$, $T>0$, is a {\em mild solution} of (\ref{07022014-1007}) if, for any $t\in [0,T)$,
\begin{equation}
u(t) = U(t,0)\bar u+\int_{0}^{t} U(t,s) F(s,u(s))\d s, \nonumber
\end{equation}
where $\{U(t,s) \}_{t\geq s\geq 0}$ is the evolution system determined
by the family $\{A(t) \}_{t\geq 0}$.  Due to standard techniques in
theory of abstract evolution equations (see \cite{Cholewa} or
\cite{Henry}), under the above assumptions, the problem
(\ref{07022014-1007}) admits a unique global mild solution $u\in C([0,+\infty),X^\alpha)$. In this section we shall deal with continuous dependence of solutions on $A(t)$, $F$ and initial values.
\begin{remark}\label{16102016-1640} {\em
Assume that $u:[0,T]\to X^\alpha$ is a mild solution of (\ref{07022014-1007}) with $T>0$ and suppose that there is $K>0$ such that for any and $s,t \geq 0$ with $t>s$
\begin{align}
\|U(t,s)\|_{\mathcal{L}(X^\alpha, X^\alpha)} \leq K   \ \ \mbox{ and } \ \
\|U(t,s)\|_{{\mathcal{L}}(X,X^\alpha)} \leq  K(1+(t-s)^{-\alpha}).\nonumber
\end{align}
Then clearly, by $(\mathcal{F}_1)$, there is a constant $\tilde C=\tilde C(\alpha, K, C ,T)>0$ such that for all $t\in (0,T]$
\begin{eqnarray*}
\|u(t)\|_\alpha & \leq & K \|\bar u\|_\alpha+\int_{0}^{t} K (1+ (t-s)^{-\alpha})\|F(s,u(s))\|\d s\\
& \leq & \tilde C(1+\|\bar u\|_\alpha) + \tilde C \int_{0}^{t}(t-s)^{-\alpha} \|u(s)\|_\alpha\d s.
\end{eqnarray*}
This in view of \cite[Lemma 1.2.9]{Cholewa} implies that there exists $\bar C=\bar C(\alpha, K, C ,T)>0$ such that
\begin{equation}
\|u(t)\|_\alpha \leq \bar C(1+\|\bar u\|_\alpha) \ \mbox{ for all }\  t\in [0,T].
\end{equation}
for all $t\in [0,T]$. \hfill $\square$
}\end{remark}

From now on we assume that families of operators $\{A_n (t)\}_{t\geq 0}$, $n\geq 0$,
satisfy conditions $(\mathcal{P}_0)-(\mathcal{P}_3)$ with commmon
spaces $X^1, X^\alpha$, constants $C_{1,\alpha}, C_{2,\alpha}$,
$C_1$, $C_2$, $M$ (not necessarily $L$) and denote by 
$\{U_n (t,s)\}_{t\geq s\geq 0}$ the evolution system determined by $\{A_n(t)\}_{t\geq 0}$. 

\begin{theorem}  \label{28012014-2357}
Under the above assumptions assume additionally that
\begin{enumerate}
	\item[$(\overline{\cal P}_1)$] there exists $K>0$ such that, for any $s,t \geq 0$ with $t>s$ and $n\geq 1$,
	\begin{equation}
	\|U_n(t,s)\|_{{\cal L}(X^\alpha,X^\alpha)} \leq K \ \ \mbox{ and } \ \ 
	\|U_n(t,s)\|_{{\cal L}(X,X^\alpha)}\leq K(1+(t-s)^{-\alpha}); \nonumber
	\end{equation}
	\item[$(\overline{\cal P}_2)$] 
	for any $\bar u\in X^\alpha$, $U_n (t,s)\bar u \to U_0 (t,s)\bar u$ in $X^\alpha$
	as $n\to +\infty$, uniformly with respect to $t \in [0,T]$ and $s\in[0,t]$,  for any fixed $T>0$.
\end{enumerate} 
Let $F_n: [0,T] \times X^\alpha \to X$, $n \geq 0$, be continuous mappings satisfying $(\mathcal{F}_1)$ and $(\mathcal{F}_2)$ with common constants $Q$ and $C$ and additionally suppose that
\begin{enumerate}
\item[$(\mathcal{F}_3)$] for each $\bar u\in X^\alpha$
 $$\int_{0}^{t} U_0 (t,s) F_n (s,\bar u) \d s \to \int_{0}^{t} U_0 (t,s) F_0(s,\bar u) \d s  \ \ \mbox{ in } \ \ X^\alpha, \mbox{ as } n \to +\infty, $$ uniformly with respect to $t \in [0,T]$;
\item[$(\mathcal{F}_4)$] the set $\{F_n(t, u)\  | \  t\in [0,T], \ n\geq 1\}$ is relatively compact in $X$ for any $u \in X^\alpha$.
\end{enumerate}
Let $(\bar u_n)$ in  $X^\alpha$ be bounded and $u_n:[0,T]\to X^\alpha$, $n\geq 1$, be mild solutions of 
\begin{equation*}  
\left\{   \begin{array}{l}
\dot{u}(t)= -A_n (t) u(t) + F_n (t, u(t)), \qquad t >  0,\\
u(0)= \bar u_n.
\end{array}  \right.
\end{equation*}
Then
\begin{enumerate}
\item[{\em (i)}] if  $\bar u_n \to \bar u_0$ in $X$, as $n \to +\infty$, for some $\bar u_0 \in X^\alpha$ then
$u_n(t) \to u_0(t)$ in $X^\alpha$, as $n \to +\infty$, uniformly with
respect to $t$ from compact subsets of $(0,T]$, where $u_0 : [0,T) \to X^\alpha$ is 
a mild solution 
\begin{equation*}  
\left\{   \begin{array}{l}
\dot{u}(t)= -A_0 (t) u(t) + F_0 (t, u(t)), \qquad t >  0,\\
u(0)= \bar u_0;
\end{array}  \right.
\end{equation*}
\item[{\em (ii)}] if $(\bar u_n)$  in $X^\alpha$ is such that $\bar u_n \to \bar u_0$ in $X^\alpha$, as $n \to +\infty$, for some $\bar u_0 \in X^\alpha$, then $u_n(t) \to u_0(t)$ in $X^\alpha$ for all $t\in [0,T]$, as $n \to +\infty$ where $u_0$ is as above.
\end{enumerate}
\end{theorem}
\begin{lemma}\label{18102016-1754}
Under the assumptions of Theorem  \ref{28012014-2357}, suppose that ${\mathcal F} \subset X$ is relatively compact and $(f_n)$ in $L^\infty ([0,T],X)$ is such that $f_n(t)\in {\mathcal F}$ for a.e. $t\in [0,T]$. Then
$$
\int_{0}^{t} \left[U_n(t,s)-U_0(t,s) \right]f_n(s)\d s \to 0 \ \mbox{ in } \ X^\alpha
$$
uniformly with respect to $t\in [0,T]$.
\end{lemma}
\noindent {\bf Proof:} First observe that, by  $(\overline{\mathcal{P}}_1)$,
$$
\|U_n(t,s)-U_0(t,s)\|_{{\cal L}(X,X^\alpha)} \leq (K+M) \big(1+(t-s)^{-\alpha}\big) \leq 2K\big(1+\delta^{-\alpha}) \mbox{ if } t > \delta, 0 \leq s \leq t-\delta,
$$
where $\delta>0$. Using this estimation one can show that $(\overline{\mathcal{P}}_2)$ implies that, for any $\eta>0$ and $\delta>0$, there exists $N$ such that, for all $n\geq N$,  one has
\begin{align}\label{18102016-1826}
\|U_n(t,s)\bar u - U_0(t,s)\bar u\|_\alpha \leq \eta  \ \mbox{ for }
 \ \bar u \in {\mathcal F}, \ t\in [\delta,T], \ s\in [0,t-\delta].
\end{align}
Now take any $\varepsilon>0$ and apply the above remark  to $\eta:= \varepsilon/2T$
and $\delta>0$ such that
$$
2KR_{\cal F}(\delta+\delta^{1-\alpha}/(1-\alpha)) <\varepsilon/2.
$$
where $R_{\mathcal F}>0$ is such that ${\mathcal F}\subset B(0,R_{\mathcal F})$. If $t\in [0,\delta]$, then clearly
\begin{align}
\left\|\int_{0}^{t}\left[U_n(t,s)-U_0(t,s)\right]f_n(s)\d s \right\|_\alpha
&\leq 2KR_{\mathcal F}\int_{0}^{t} \big(1+(t-s)^{-\alpha}\big)\d s \nonumber \\
&\leq 2KR_{\mathcal F}\big(\delta+ \delta^{1-\alpha}/(1-\alpha)\big)< \varepsilon/2  < \varepsilon.\nonumber
\end{align}
If $t\in [\delta, T]$, then
\begin{align}\label{18102016-1813}
\left\|\int_{t-\delta}^{t} \left[U_n(t,s)-U_0(t,s)\right]f_n(s)\d s \right\|_\alpha &\leq 2KR_{\mathcal F}\int_{t-\delta}^{t}\big(1+ (t-s)^{-\alpha}\big)\d s \nonumber \\
&=  2KR_{\mathcal F} \big(\delta+ \delta^{1-\alpha}/(1-\alpha)\big) < \varepsilon/2.
\end{align}
and, for all $n\geq N$,
$$
\left\|\int_{0}^{t-\delta} \left[U_n(t,s)-U_0(t,s)\right]f_n(s)\d s\right\|_\alpha
\leq \eta (t-\delta)= \frac{\varepsilon}{2}\cdot \frac{t-\delta}{T}  < \varepsilon/2,
$$
which, together with \eqref{18102016-1813}, yields
$$
\left\|\int_{0}^{t} \left[U_n(t,s)-U_0(t,s)\right]f_n(s)\d s\right\|_\alpha <\varepsilon
$$
for all $n\geq N$. \hfill $\square$
\\ \ \\
\noindent{\bf Proof of Theorem \ref{28012014-2357}:}
Observe that
\begin{align}
u_n(t)-u_0(t)= U_n(t,0)\bar u_n - U_0(t,0) \bar u_0 & + \int_0^t
U_n(t,\tau)F_n(\tau, u_n(\tau)) \d\tau \nonumber \\
& - \int_0^t U_0(t,\tau)F_0(\tau,u_0(\tau)) \d\tau. \nonumber
\end{align}
Clearly, by use of $(\overline{\mathcal{P}}_1)$,
\begin{align}
\|U_n(t,0)\bar u_n - U_0(t,0) \bar u_0\|_\alpha & \leq \|U_n(t,0)
(\bar u_n - \bar u_0) \|_\alpha + \| U_n(t,0)\bar u_0 - U_0(t,0)
\bar u_0 \|_\alpha \nonumber  \\
& \leq K(1+ t^{-\alpha}) \|\bar u_n - \bar
u_0 \| + \| U_n(t,0)\bar u_0 - U_0(t,0) \bar u_0 \|_\alpha.
\nonumber
\end{align}
In consequence, if $\bar u_n \to \bar u_0$ in $X$, then $(\overline{\mathcal{P}}_2)$ implies that
\begin{equation}\label{04022014-0035}
U_n(t,0)\bar u_n \to U_0(t,0) \bar u_0\   \mbox { in }  \ X^\alpha,
\mbox{ uniformly on compact subsets of } (0,T].
\end{equation}
Further, note that
\begin{equation}
\int_0^t \bigg(
U_n(t,\tau)F_n(\tau,u_n(\tau))-U_0(t,\tau)F_0(\tau,u_0(\tau))\bigg)\d\tau
= I_{1,n}(t)+ I_{2,n}(t)+ I_{3,n}(t) \nonumber
\end{equation}
where
\begin{align}
I_{1,n}(t) &:= \int_0^t \bigg(U_n(t,\tau)F_n(\tau,u_n(\tau))-U_n(t,\tau)F_n(\tau,u_0(\tau))\bigg)
\d\tau, \nonumber \\
I_{2,n}(t) &:= \int_0^t \bigg(
U_n(t,\tau)F_n(\tau,u_0(\tau))-U_0(t,\tau)F_n(\tau,u_0(\tau))\bigg)\d\tau, \nonumber \\
I_{3,n}(t) &:= \int_0^t \bigg(U_0(t,\tau)F_n(\tau,u_0(\tau))-U_0(t,\tau)F_0(\tau,u_0(\tau))\bigg)\d\tau.\nonumber
\end{align}
\noindent  In view of $(\overline{\mathcal{P}}_1)$ and $(\mathcal{F}_1)$, for any $t \in [0,T]$,
\begin{align}
\| I_{1,n} (t)\|_\alpha &\leq \int_0^t \| U_n(t,\tau)
\|_{\mathcal{L}(X,X^\alpha)} \| F_n (\tau,u_n(\tau))- F_n(\tau,
u_0(\tau) )\| \d\tau \nonumber \\
& \leq KQ \int_0^t
\big(1+(t-\tau)^{-\alpha}\big)\| u_n (\tau)-u_0(\tau)
\|_\alpha \d\tau \nonumber \\
& \leq KQ(T^\alpha+1) \int_0^t (t-\tau)^{-\alpha}\| u_n (\tau)-u_0(\tau)
\|_\alpha \d\tau. \nonumber
\end{align}
In order to estimate $I_{2,n}(t)$, we need to know that the set
$\mathcal{F}:= \{F_n (\tau,u_0 (\tau))\  | \ \tau \in [0,T],\  n \geq
1 \}$
is relatively compact in $X$. To see this, take any sequence
$(v_k)$ in $\mathcal{F}$. Then, for each $k \geq 1$ there exist $n_k\in \N$\  and $\tau_k \in [0,T]$ such that $v_k = F_{n_k}(\tau_k, u_0(\tau_k))$.
Without loss of generality  we may assume that $\tau_k \to \tau_0$, as $k \to +\infty$ for some $\tau_0 \in [0,T]$.
By  continuity, $u_0 (\tau_k) \to u_0 (\tau_0)$ in $X^\alpha$.
Since the functions $F_n$ have common constants $C$ and $Q$ in $({\cal F}_1)$ and $({\cal F}_2)$, in view of $(\overline{\mathcal{P}}_1)$, Remark \ref{16102016-1640} and
the boundedness of $(\bar u_n)$ in $X^\alpha$,  for any $\varepsilon >0$ one can find $k_\varepsilon >0$ such that
$$
\| F_{n_k}(\tau_k, u_0(\tau_k))-F_{n_k}(\tau_k, u_0 (\tau_0))\| < \varepsilon \ \mbox{ for all } k\geq k_\varepsilon,
$$
that is
$$
\{v_k \ | \ k\geq k_\varepsilon \} \subset
\{F_n(\tau,u_0(\tau_0))
\ | \  \tau \in [0,T],\ n\geq 1 \} + B_X (0,\varepsilon).
$$
Hence, in view of $(\mathcal{F}_4)$, $(v_k)$ contains
a convergent subsequence, which proves that  $\mathcal{F}$ is relatively
compact in $X$. This according to Lemma \ref{18102016-1754} implies that $\|I_{2,n} (t)\|_\alpha \to 0$, as $n \to +\infty$, uniformly with respect to $t \in [0,T]$.\\
\indent Moreover, one may prove by use of $(\overline{\mathcal{P}}_3)$ that $\|I_{3,n}(t)\|_\alpha\to 0$, as $n \to +\infty$, uniformly with respect to $t\in [0,T]$.\\
\indent Summing up, we have
\begin{equation}
\|u_n(t) - u_0(t) \|_\alpha \leq \gamma_n(t) +  KQ(T^\alpha+1) \int_0^t (t-\tau)^{-\alpha}\| u_n (\tau)-u_0(\tau)
\|_\alpha \d\tau \nonumber
\end{equation}
where $\gamma_n(t):= \|U_n(t,0)\bar u_n - U_0(t,0) \bar u_0\|_\alpha
+ \| I_{2,n}(t)\|_\alpha + \| I_{3,n}(t)\|_\alpha$. Note that $\gamma_n(t) \to 0$, as $n \to +\infty$,  uniformly with respect to $t$ from compact subsets of $(0,T]$. Furthermore, by use of Lemma 7.1.1 of \cite{Henry}, we get
\begin{equation}\label{18102016-2152}
\|u_n(t)-u_0(t)\|_\alpha \leq \gamma_n(t) + \bar K \int_{0}^{t}
(t-\tau)^{-\alpha}\gamma_n(\tau)\d \tau
\end{equation}
for some constant $\bar K>0$.
Taking an arbitrary $\delta\in (0,T]$, we get, for $t\in [\delta, T]$,
\begin{align*}
\int_{0}^{t} (t-\tau)^{-\alpha}\gamma_n(\tau)\d \tau \leq
\frac{2^\alpha}{\delta^\alpha} \int_{0}^{t-\delta/2} \gamma_n(\tau)\d \tau +
\int_{t-\delta/2}^{t} (t-\tau)^{-\alpha}\gamma_n(\tau)\d \tau\\
\leq
\frac{2^\alpha}{\delta^\alpha} \int_{0}^{T} \gamma_n(\tau)\d \tau
+ \frac{T^{1-\alpha}}{1-\alpha}\cdot
\max_{s\in [\delta/2,T]} \gamma_n (s).
\end{align*}
Since $\gamma_n(t)\leq \bar M(1+t^{-\alpha})$, for $t\in [0,T]$,
with some constant $\bar M > 0$, we infer, by the dominated convergence
theorem, that $\int_{0}^{T}\gamma_n (\tau)\d \tau \to 0$, as $n\to+\infty$. This together with the uniform convergence of $\gamma_n$ on $[\delta/2,T]$ proves that
$$
\int_{0}^{t} (t-\tau)^{-\alpha}\gamma_n(\tau)\d \tau \to 0
$$
uniformly with respect to $t\in [\delta,T]$. Finally, we get $\|u_n(t)-u_0(t)\|_\alpha\to 0$, as $n \to +\infty$, uniformly with respect to $t\in [\delta,T]$.\\
\indent In case (ii), $(\overline{\mathcal{P}}_1)$ implies that
\begin{align}
\|U_n(t,0)\bar u_n - U_0(t,0) \bar u_0\|_\alpha  \leq K \|\bar u_n -
\bar u_0 \|_\alpha + \| U_n(t,0)\bar u_0 - U_0(t,0) \bar u_0
\|_\alpha, \nonumber
\end{align}
which means that $\|U_n(t,0)\bar u_n - U_0(t,0) \bar u_0\|_\alpha\to 0$, as $n \to +\infty$, uniformly with respect to $t\in [0,T]$. Hence $\gamma_n$ converges to $0$
uniformly on $[0,T]$. This together with \eqref{18102016-2152} leads to (ii).\hfill $\square$

\section{Averaging principle for nonlinear problems}
In this section we shall apply the continuity results obtained in the previous one. 
Consider, like in Section 1, the family of linear operators $\{A(t)\}_{t \geq 0}$ in a Banach space $X$ satisfying conditions
$(\mathcal{P}_0)$--$(\mathcal{P}_3)$, ($\mathcal{\widehat{P}}_1)$ and ($\mathcal{\widehat{P}}_2)$. For the future reference we shall deal with the case 
where $F$ depends on parameter that is we assume that $F: [0, +\infty) \times X^\alpha \times P \to X$, where $P$ is a metric space of parameters, is a continuous
mapping such that $F(\cdot,\cdot, \mu)$, $\mu\in P$, satisfy $(\mathcal{F}_1)$ and $(\mathcal{F}_2)$ with common constants $Q$ and $C$ and a parametrized version of $(\widehat{\mathcal{F}})$ holds:
\begin{enumerate}
\item[$(\widehat{\mathcal{F}}_\mu)$]
the set
$\{F(t,\bar{u},\mu)\ |\ t \geq 0, \mu\in P\}$ is relatively compact for any    $\bar{u} \in X^\alpha$, and there is a continuous function $\widehat{F}: X^\alpha \times P \to X$ such that, for any $\bar{u} \in X^\alpha$ and $\mu \in P$,
\begin{equation}
\lim_{\omega \to +\infty,\ \nu \to \mu} \frac{1}{\omega}
\int_0^\omega F(\tau, \bar{u}, \nu) \d\tau = \widehat{F}(\bar{u}, \mu)\ \
\text{ in } X. \nonumber
\end{equation}
\end{enumerate}
Define $\{A^{(\lambda)}(t) \}_{t \geq 0}$,
$\lambda > 0$, as in $(\ref{20012017-2206})$ and let $F^{(\mu,
\lambda)}:[0,+\infty) \times X^\alpha \to X$,\  $\mu \in P$,
$\lambda
>0$, be given by
\begin{displaymath}
F^{(\mu, \lambda)}(t, \bar{u}):=  F(t/ \lambda, \bar{u}, \mu),\ t
\geq 0, \bar{u} \in X^\alpha.
\end{displaymath}

\noindent Consider the perturbed problem
\begin{align}              \label{05022014-0054}
\left\{
\begin{array}{ll}
\dot{u}(t)= -A^{(\lambda)}(t)u(t)+ F^{(\mu, \lambda)}(t,u(t)), & t>0, \\
u(0)= \bar{u}, &
\end{array}
\right.
\end{align}
where $\bar u \in X^\alpha$, $\mu \in P$ and $\lma>0$.
\noindent By $u(\cdot;\bar{u}, \mu, \lambda): [0, +\infty) \to
X^\alpha$ we denote its mild solution. The following parametrized averaging principle for nonlinear parabolic problems (a version of Theorem \ref{20012017-2210}) holds.
\begin{theorem} \label{06012017-2015}
Under the above assumptions, if $(\mu_n)$ in $P$ and $(\lma_n)$ in $(0,+\infty)$ are
sequences such that $\mu_n \to \mu_0$, $\lma_n \to
0^+$ as $n \to +\infty$ for some $\mu_0
\in P$, $(\bar u_n)$ is a bounded sequence in $X^\alpha$ and $u_n : [0, +\infty) \to X^\alpha$, $n\geq 1$, are solutions of
$(\ref{05022014-0054})$ with $\lma:=\lma_n$, $\mu:=\mu_n$ and $\bar
u:= \bar u_n$, respectively, then
\begin{enumerate}
\item[{\em (i)}]
if $\bar{u}_n \to \bar{u}$ in $X$ then
\begin{equation}
u(t; \bar{u}_n, \mu_n, \lma_n) \to \widehat{u}(t; \bar{u}_0, \mu_0)
\mbox{ in } X^\alpha, \mbox{as } n \to +\infty, \nonumber
\end{equation}
uniformly for $t$ from  compact subsets of $[0,+\infty)$, where $\widehat u (\cdot;\bar u_0,\mu_0) : [0, +\infty) \to X^\alpha$ is the solution
of
\begin{align}
\left\{
\begin{array}{ll}
\dot{u}(t)= - \widehat{A}u(t)+ \widehat F(u(t), \mu_0), \ t>0,\\
u(0)= \bar{u}_0; &
\end{array}\right.
\end{align}
\item[{\em (ii)}] if $\bar{u}_n \to \bar{u}$ in $X^\alpha$ then
\begin{equation}
u(t; \bar{u}_n, \mu_n, \lma_n) \to \widehat{u}(t; \bar{u}_0, \mu_0)
\mbox{ in } X^\alpha \mbox{as } n \to +\infty, \nonumber
\end{equation}
\noindent uniformly for $t$ from compact subsets of $(0,+\infty)$.
\end{enumerate}
\end{theorem}
\begin{lemma}[see Lemma 3.4 in \cite{Cwiszewski-RL}]   \label{26102016-1116}
Let $T>0$ and $F_n:[0,T]\times X^\alpha \to X$, $n \geq 0$ satisfy $(\mathcal{F}_1)$ and $(\mathcal{F}_2)$ with common constants $Q$ and $C$ (independent of $n$) and let for each $\bar u \in X^\alpha$
\begin{equation}
\int_0^t F_n(s,\bar u)\d s \to \int_0^t F_0(s,\bar u) \d s \textnormal{ in }X^\alpha \textnormal{ as } n\to +\infty, \nonumber
\end{equation}
uniformly with respect to $t \in [0,T]$. Then, for any continuous  $u:[0,T]\to X^\alpha$
$$
\int_{0}^{t} \widehat S(t-s) F_n (s, u(s)) \d s \to \int_{0}^{t} \widehat S(t-s) F_0 (s, u(s)) \d s \ \mbox{ in } X^\alpha, \mbox{ as } \ n
\to +\infty,
$$
uniformly with respect to $t\in [0,T]$.\hfill $\square$
\end{lemma}

\noindent {\bf Proof of Theorem \ref{06012017-2015}:} We shall apply Theorem \ref{28012014-2357}. Let $\{U^{(\lambda)}(t,s)\}_{t\geq s\geq 0}$ be the evolution system generated by $\{A^{(\lambda)}(t)\}_{t\geq 0}$,
$\lambda>0$. Put $A_n(t) := A^{(\lambda_n)}(t)$ and
$U_n(t,s) := U^{(\lambda_n)}(t,s)$, $n\geq 1$. It is clear that ($\mathcal{\widehat{P}}_1)$ implies $(\overline{\mathcal{P}}_1)$ of Theorem \ref{28012014-2357}. Further, if we fix $T>0$ and then, due to Theorem \ref{17012017-1853}, for any $\bar u \in X^\alpha$,
$$
\lim_{n\to \infty} U_n (t,s)\bar u = \widehat S(t-s)\bar u \ \mbox{ in } X^\alpha,
$$
uniformly with respect to $t\in [0, T]$ and $s\in [0,t]$, i.e. $(\overline{\mathcal{P}}_2)$ of Theorem \ref{28012014-2357} is satisfied. Let $F_n:=F(\cdot/\lma_n,\cdot, \mu_n)$ and $F_0:=\widehat F(\cdot, \mu_0)$.  Observe that, for any $\bar u \in
X^\alpha$ and $t>0$, we get by $(\widehat{\mathcal{F}}_\mu)$,
$$
\int_{0}^{t} F_n (s,\bar u)\d s=  \lma_n \int_{0}^{t/\lma_n}
F(\rho,\bar u,\mu_n) \d \rho\to t\widehat F(\bar u,\mu_0)  \
\mbox{ in } X, \mbox{ as }  n\to+\infty,
$$
uniformly with respect to $t\in [0,T]$. This, in view of Lemma \ref{26102016-1116}, implies
$$
\int_{0}^{t}\widehat S(t-s)F_n(s,\bar u)\d s \to \int_{0}^{t}
\widehat S(t-s)\widehat F(\bar u, \mu_0)\d s, \mbox{ as }  n\to+\infty,
$$
uniformly with respect to $t\in [0,T]$, i.e. condition $(\mathcal{F}_3)$ of Theorem \ref{28012014-2357} holds. Observe also that $\{F_n(t,\bar u)\mid t\in [0,T], n\geq 1\}$ is contained in
$F([0,+\infty)\times \{\bar u\}\times P)$ and, in consequence, is compact due to  $(\widehat{\mathcal{F}}_\mu)$, i.e. $(\mathcal{F}_4)$ of Theorem \ref{28012014-2357} is also satisfied. Finally, we get both the assertions by a direct use of Theorem \ref{28012014-2357}. \hfill $\square$

\section{Applications to parabolic PDE's}
Consider the equation
\begin{align} \label{07012017-0031}
\left\{
\begin{array}{ll}
\displaystyle{\frac{\partial u}{\partial t}(x,t)}= \displaystyle{\sum_{i,j=1}^N a_{ij}(t/\lma) \frac{\partial^2 u}{\partial x_j \partial x_i} (x,t)}- g(t/\lma)u(x,t)+ f(t/\lma,x),\ x \in \mathbb{R}^N, t>0,\\
u(x,0)= \bar{v}(x)
\end{array}
\right.
\end{align}
\noindent where $N \geq 1$, $\lma>0$, $\bar v \in H^1(\mathbb{R}^N)$, functions $a_{i,j}: [0,+\infty) \to
\mathbb{R}$ are bounded and such that $a_{ij}=a_{ji}$ for any $t \geq 0$ and $i,j=1,\ldots N$, $g:[0,+\infty) \to \mathbb{R}$ is such that $0<\bar c_1\leq g(t)\leq \bar c_2$ for $t\geq 0$ and some $\bar c_1,\bar c_2>0$ and the following conditions hold
\begin{enumerate}
	\item[$(\mathcal{Z}_1)$] there exists a positive constant $\nu_0>0$ such that
	\begin{align}
	\sum_{i,j=1}^N a_{ij}(t) \xi_i \xi_j \geq \nu_0|\xi|^2 \textnormal{ for all }t \geq 0 \textnormal{ and } \xi \in \mathbb{R}^N;\nonumber
	\end{align}
	\item[$(\mathcal{Z}_2)$] there are constants $\bar{L}>0$ and $0 \leq \gamma <1$ such that for any $s, t\geq 0$ and $i,j=1,\ldots,N$
	\begin{equation}
	\left| a_{ij}(t)-
	a_{ij}(s) \right| \leq \bar{L} \left| t-s \right|^{\gamma} \textnormal{ and } \left| g(t)-g(s)\right| \leq \bar{L} \left| t-s \right|^{\gamma};
	\nonumber
	\end{equation}
	\item[$(\mathcal{Z}_3)$] there are $\widehat a_{ij} \in \mathbb{R}$, $i,j=1,\ldots,N$, and $\widehat g \in \mathbb{R}$ such that
	\begin{equation}
	\lim_{\omega \to +\infty} \frac{1}{\omega} \int_0^\omega a_{ij}(t) \d t = \widehat a_{ij},\ \ \ \ \lim_{\omega \to +\infty} \frac{1}{\omega} \int_0^\omega g(t) \d t = \widehat g. \nonumber
	\end{equation}
\end{enumerate}
We assume that $f:[0,+\infty) \times \mathbb{R}^N \to \mathbb{R}$ is continuous and such that
\begin{equation*}
\left|f(t,x)\right| \leq k(x) \textnormal{ for all } x \in \mathbb{R}^N,\ t\geq 0, \leqno{({\cal Z}_4)}
\end{equation*}
for some $k \in L^2(\mathbb{R}^N)$,
\begin{equation*}
\int_{\mathbb{R}^N} |f(t,x+y)-f(t,x)|^2 \d x \to 0 \textnormal{ as }y\to 0, \textnormal{ uniformly in } t \geq 0 \leqno{({\cal Z}_5)}
\end{equation*}
and there is $\widehat f \in L^2(\mathbb{R}^N)$ such that
\begin{equation}
\ \ \ \ \ \ \ \ \ \ \ \ \ \ \ \widehat f(x) = \lim_{\omega \to +\infty} \frac{1}{\omega}\int_0^\omega f(t,x) \d t \textnormal{ for a.e. }x \in \mathbb{R}^N.\ \ \ \ \ \ \ \ \ \ \ \ \ \ \ \ \ \ \ \ \ \  \ \ \ \ \ \ \ \ \ \ \ \ \ \nonumber \leqno{({\cal Z}_6)}
\end{equation}
\indent
Define the family of linear operators $\{A(t): D \to X\}_{t \geq 0}$, where
$X:=L^2(\mathbb{R}^N)$ and $D:=H^2(\mathbb{R}^N)$, by
\begin{equation}
A(t)u:= -\sum_{i,j=1}^N a_{ij}(t) \frac{\partial^2 u}{\partial x_j \partial x_i}+ g(t)u \ \textnormal{ for } \ u\in D. \nonumber
\end{equation}
We shall first demonstrate that the family $\{A(t)\}_{t \geq 0}$ satisfies conditions $(\mathcal{P}_1)$--$(\mathcal{P}_3)$. The computations are rather standard and the main issue here is the independence of the constants of $t$ and $\lma$.

\begin{lemma}\label{07022017-1248}
For each $t \geq 0$, $\sigma(A(t)) \subset \{\zeta \in \C \mid \mathrm{Re}\, \zeta > 0 \}$ and there is $M>0$ such that for all $t \geq 0$ and $\zeta \in \C$ with $\mathrm{Re}\ \zeta \leq 0$
	\begin{equation}
	\|(\zeta I - A(t))^{-1}\|_{\mathcal{L}(L^2, L^2)} \leq \frac{M}{1+\left|\zeta\right|}.\nonumber
	\end{equation}
\end{lemma}
\noindent {\bf Proof:} 
Let us fix $t\geq 0$ and take an arbitrary $h \in L^2(\mathbb{R}^N)$. For $\zeta \in \C$ with $\mathrm{Re}\ \zeta \leq 0$ and $u \in H^2(\mathbb{R}^N)$ consider
\begin{equation}
-A(t)u + \zeta u = h,\ t\geq 0.\nonumber
\end{equation}
Then, by use of the Fourier transform\footnote[1]{Recall that, for  $v \in L^1(\mathbb{R}^N)\cap L^2(\mathbb{R}^N)$,
	$$
	(\mathcal{F}v)(\xi)=\frac{1}{(2\pi)^{N/2}}\int_{\mathbb{R}^N}e^{-ix\cdot \xi}v(x)\d x.
	$$} 
one can get 
\begin{equation}
\left(-\sum_{i,j=1}^N a_{ij}(t) \xi_i \xi_j -g(t) + \zeta \right)(\mathcal{F}u)(\xi) = (\mathcal{F}h)(\xi), \ \xi\in\R^N,\nonumber
\end{equation}
where $\mathcal{F}v$ stands for the Fourier transform of $v \in L^2(\mathbb{R}^N)$. This implies that
\begin{equation}
\int_{\mathbb{R}^N} \left(\sum_{i,j=1}^N a_{ij}(t) \xi_i \xi_j +g(t) - \zeta \right)\left|(\mathcal{F}u)(\xi)\right|^2 \d \xi = \big(\mathcal{F}h,\mathcal{F}u\big)_{L^2}\nonumber
\end{equation}
Hence, by use of $({\cal Z}_1)$ and the Cauchy-Schwarz inequality,  we get
\begin{equation}
|\bar c_1 -\zeta|\cdot \|u\|^2_{L^2} \leq \|h\|_{L^2}\|u\|_{L^2}\nonumber
\end{equation}
and  in consequence\footnote{Observe that, for any positive real number $a$ and complex number $b$ with $\mathrm{Re\ } b \leq 0$, one has
	$$
	\left| a-b\right| \geq (1/\sqrt{2})(a+\left| b \right|).
	$$} 
\begin{equation}
\frac{\left(\bar c_1 +|\zeta| \right)}{\sqrt{2}} \|u\|_{L^2} \leq \|h\|_{L^2},\nonumber
\end{equation}
which ends the proof.\hfill $\square$\\

\begin{lemma}\label{07022017-0704}
There are $C_1, C_2>0$ such that
$$
C_1\|u\|_{H^2} \leq \|A(t)u\|_{L^2} \leq C_2 \|u\|_{H^2} \text{ for all } u\in H^2(\R^N) \text{ and  } t\geq 0,
$$
i.e. $(\mathcal{P}_1)$ with $X^1:=D=H^2(\R^N)$ is satisfied by the operators $A(t)$, $t \geq 0$.
\end{lemma}
\noindent {\bf Proof:} Observe that for each $t \geq 0$ and $u \in H^2(\mathbb{R}^N)$
\begin{align}
\left\|A(t)u\right\|_{L^2} &\leq \left\|\sum_{i,j=1}^N a_{ij}(t)\frac{\partial^2 u}{\partial x_j \partial x_i}\right\|_{L^2} + \left\|g(t)u\right\|_{L^2} \leq 
a_\infty \sum_{i,j=1}^N \left\| \frac{\partial^2 u}{\partial x_j \partial x_i}\right\|_{L^2} + \bar c_2 \|u\|_{L^2} \nonumber \\ 
&\leq \max\{a_\infty, \bar c_2\} \|u\|_{H^2}, \label{21012017-1945}
\end{align}
where $a_\infty>0$ is such that $|a_{ij}(t)|\leq a_\infty$ for $t \geq 0$ and $i,j=1,\ldots,N$. On the other hand, for $t \geq 0$ and $u \in H^2(\mathbb{R}^N)$ one has 
\begin{align}
\|A(t)u\|^2_{L^2} & = \left\|\mathcal{F}A(t)u\right\|^2_{L^2} = \int_{\mathbb{R}^N}|\sum_{i,j=1}^N a_{ij}(t)\xi_i \xi_j+g(t)|^2\left|(\mathcal{F}u)(\xi)\right|^2 \d \xi \nonumber \\
& \geq \int_{\mathbb{R}^N}\left(\nu_0 |\xi|^2 +\bar c_1 \right)^2 \left|(\mathcal{F}u)(\xi)\right|^2 \d \xi \geq
\int_{\mathbb{R}^N}\left(\nu_0^2 |\xi|^4 +\bar c_1^2 \right) \left|(\mathcal{F}u)(\xi)\right|^2 \d \xi \nonumber \\
& \geq \beta \|u\|^2_{H^2}\nonumber
\end{align} 
for some $\beta >0$, which yields the assertion. \hfill $\square$\\

\begin{lemma}
There are constants $\bar L>0$ and $\gamma \in (0,1]$ such that for $A(t)$, $t\geq 0$, satisfy $(\mathcal{P}_3)$.
\end{lemma}
\noindent {\bf Proof:} Take any $w \in L^2(\mathbb{R}^N)$ and let  $u= A(r)^{-1}w$. Then $u \in H^2(\mathbb{R}^N)$ and by $(\mathcal{Z}_1)$, for any $t,s,r \geq 0$ we have 
\begin{align}
\|&(A(t)-A(s))A(r)^{-1}w\|_{L^2} = \|A(t)u-A(s)u\|_{L^2} \nonumber \\ & \leq \left(\int_{\mathbb{R}^N} \left|\sum_{i,j=1}^N (a_{ij}(t)-a_{ij}(s)) \frac{\partial^2 u}{\partial x_j \partial x_i} \right|^2 \d x \right)^{1/2} +\left( \int_{\mathbb{R}^N}\left| g(t)-g(s) \right|^2 \left| u \right|^2 \d x \right)^{1/2}\nonumber \\ &
\leq  \bar L |t-s|^{\gamma} \|u\|_{H^2} =  \bar L |t-s|^{\gamma} \|A(r)^{-1}w\|_{H^2}, \nonumber
\end{align}
which gives
\begin{equation}
\|(A(t)-A(s))A(r)^{-1}w\|_{L^2} \leq \bar L |t-s|^{\gamma} \|A(r)^{-1}w\|_{H^2} \leq 
(\bar L/C_1)|t-s|^{\gamma} \| w\|_{L^2}
\nonumber
\end{equation}
where $C_1>0$ is the constant from Lemma \ref{07022017-0704}. \hfill $\square$\\

In consequence, $A(t)$ are sectorial operators and, since conditions $({\cal P}_1)-({\cal P}_3)$ hold, the parabolic equation
\begin{equation}
\dot{u}(t) = -A(t) u(t),\ t>0, \nonumber
\end{equation}
generates a parabolic evolution system $\{U(t,s)\}_{t \geq s \geq 0}$ on $L^2(\mathbb{R}^N)$. It is left to explain the equivalence of the norms in the $1/2$-fractional spaces determined by $A(t)$ (for each $t \geq 0$).
\begin{lemma}
For each $t \geq 0$, $A(t)$ determines the fractional space $X^{1/2}=H^1(\mathbb{R}^N)$ and there are $C_{1,1/2}, C_{2,1/2}>0$ such that
$$
C_{1,1/2}\|u\|_{H^1} \leq \|A(t)^{1/2}u\|_{L^2} \leq C_{2,1/2}\|u\|_{H^1},
$$
i.e. $(\mathcal{P}_0)$ holds.
\end{lemma}
\noindent {\bf Proof:} Fix $t\geq 0$. Clearly by $(\mathcal{Z}_2)$ it follows that, for any $u \in H^2(\mathbb{R}^N)$,
\begin{align}
\|A(t)^{1/2}u\|^2_{L^2}&= \left| \left(A(t)u,u\right)_{L^2}\right| = \left| \int_{\mathbb{R}^N} \sum_{i,j=1}^N a_{ij}(t) \frac{\partial u}{\partial x_i}\frac{\partial u}{\partial x_j} + g(t) u^2 \d x \right| \nonumber \\
& \geq \int_{\mathbb{R}^N} \nu_0 \left| \nabla u \right|^2 + \bar c_1 \left| u \right|^2 \d x.\nonumber
\end{align}
\noindent This gives the existence of $C_{1, 1/2}>0$ such that, for all $u \in H^2(\mathbb{R}^N)$,
\begin{equation}
C_{1, 1/2}\|u\|_{H^1} \leq \|A(t)^{1/2}u\|_{L^2}.\label{21012017-2251}
\end{equation}
On the other hand it is easy to see that
\begin{align}
\|A(t)^{1/2}u\|^2_{L^2}
\leq N^2 a_\infty \|\nabla u\|^2_{L^2} + \bar c_2 \|u\|^2_{L^2} \leq C_{2,1/2}^2 \|u\|^2_{H^1},\label{21012017-2250}
\end{align}
where $C_{2,1/2}= (N^2 a_\infty+\bar c_2)^{1/2}$. Since $u\in H^2(\R^N)$ was arbitrary, (\ref{21012017-2251}) together with (\ref{21012017-2250}) imply that the $1/2$-space for $A(t)$ is equal to $H^1(\mathbb{R}^N)$.\hfill $\square$
$\mbox{ }$\\

\indent In order to apply our results,  consider the following problem, corresponding to \eqref{07012017-0031}
\begin{equation}
\dot{u}(t) = - A (t/\lambda)u(t) + F(t/\lambda), \qquad \ t>0, \label{18122016_2044}
\end{equation}
\noindent with $F:[0,+\infty) \to L^2(\mathbb{R}^N)$ given by $F(t):= f(t,
\cdot)$, for each $t \geq 0$. We claim that the solutions of the above equation converge with $\lambda \to 0^+$ to the solution of
\begin{equation}
\dot{u}(t) = -\widehat{A}u(t) + \widehat F,\ t>0, \nonumber
\end{equation}
where $\widehat A: D(\widehat A) \to L^2(\mathbb{R}^N)$ is given by
$$
\widehat A u:= -\sum_{i,j=1}^N \widehat a_{ij} \frac{\partial^2 u}{\partial x_j \partial x_i}+\widehat g u  \ \textnormal{ for } \ u\in D(\widehat A):=H^2(\mathbb{R}^N),
$$
and $\widehat F(x)=\widehat f(x)$ for a.e. $x \in \mathbb{R}^N$.
To this end we need to verify assumptions $(\widehat {\cal P}_1)$ and $(\widehat{\cal P}_2)$.\\
\indent Denoting by $\{U^{(\lambda)}(t,s)\}_{t \geq s \geq 0}$ the evolution system determined by $\{A(t/\lambda)\}_{t \geq 0}$, $\lma >0$, and  using
the Fourier transform,  we get
\begin{equation}
\mathcal{F}(U^{(\lambda)}(t,s)u)(\xi) = \exp\left(-\int_s^t \left(\sum_{i,j=1}^N a_{ij}(\rho) \xi_i \xi_j +g(\rho)\right) \d \rho\right) (\mathcal{F}u)(\xi), \ \xi\in\R^N. \nonumber
\end{equation}
It can be shown, that there is a constant $K>0$ depending only on $\nu_0$ and $\bar c_1$, such that for any $t \geq s$ and $u \in H^1(\mathbb{R}^N)$,
\begin{equation}
\|U^{(\lambda)}(t,s)u\|_{H^1} \leq K (1+(t-s)^{-1/2})\|u\|_{L^2}\label{18122016-1847}
\end{equation}
and
\begin{equation}
\|U^{(\lambda)}(t,s)u\|_{H^1} \leq K \|u\|_{H^1}\label{06012017-2352}
\end{equation}
(see \cite{Antoci-Prizzi}).
Clearly, (\ref{18122016-1847}) and (\ref{06012017-2352}) mean that $(\mathcal{\widehat{P}}_1)$ holds. Further, if $\{ \widehat{S}(t):L^2(\mathbb{R}^N) \to L^2(\mathbb{R}^N)\}_{t\geq 0}$ is an analytic $C_0$-semigroup generated by $-\widehat{A}$, then it is well known that there is a constant $\widehat M>0$ such that $\|\widehat S(t)\|_{\mathcal{L}(L^2,L^2)} \leq \widehat M$ (see e.g. \cite[Theorem1.3.2]{Cholewa}) and $(\mathcal{Z}_3)$ implies that, for any $u \in H^2(\mathbb{R}^N)$,
\begin{align}
\frac{1}{\omega}\int_0^\omega (A(t)u-\widehat A) u \d t & = \frac{1}{\omega}\int_0^\omega \left(\sum_{i,j=1}^N \left(a_{ij}(t) -\widehat a_{ij} \right)\right)\!\frac{\partial^2 u}{\partial x_j \partial x_i} \d t   \nonumber \\&+ \frac{1}{\omega}\int_0^\omega \left(g(t)-\widehat g\right) u \d t \to 0, \textnormal{ as } \omega \to +\infty,\nonumber
\end{align}
i.e. $(\mathcal{\widehat{P}}_2)$ is satisfied. 
Finally note that, due to $({\cal Z}_4)$, $({\cal Z}_5)$ and the Kolmogorov-Riesz compactness criterion, the set $F([0,+\infty))$ is relatively compact in $L^2(\mathbb{R}^N)$ and, in virtue of $({\cal Z}_6)$ one gets
\begin{equation}
\lim_{\omega \to +\infty} \frac{1}{\omega} \int_0^\omega F(t) \d t = \hat f \mbox{ in } L^2(\mathbb{R}^N). \nonumber
\end{equation}
\noindent Thus by use of Theorem \ref{20012017-2210} we have
\begin{corollary} $\mbox{ }$ 
\begin{enumerate}
\item[{\em (i)}]
If $\bar v \to \bar u$ in $L^2(\mathbb{R}^N)$, $\bar v$ come from a bounded subset of $H^1(\R^N)$ and $\bar u\in H^1(\R^N)$, then, for each $t>0$,
$u(t;\bar v) \to \widehat u(t) \textnormal{ in } H^1(\mathbb{R}^N)$,
where $u(\cdot;\bar v):[0;+\infty)\to H^1(\mathbb{R}^N)$ is the mild solution of \eqref{07012017-0031} and $\widehat u:[0;+\infty)\to H^1(\mathbb{R}^N)$ is the mild solution of 
\begin{align}
\left\{
\begin{array}{ll}
\displaystyle{\frac{\partial u}{\partial t}}(x,t)= \displaystyle{\sum_{i,j=1}^N} \widehat a_{ij} \frac{\partial^2 u}{\partial x_j \partial x_i} (x,t)-\widehat g\, u(x,t) + \widehat f(x), \ \  x \in \mathbb{R}^N,\ t>0\\\nonumber
u(x,0)= \bar{u}(x)
\end{array}
\right.
\end{align}
and the convergence is uniform with respect to $t$ from compact subsets of $(0,+\infty)$;
\item[{\em (ii)}] If $\bar v \to \bar u$ in $H^1(\mathbb{R}^N)$, then, for each $t>0$, $u(t;\bar v) \to \widehat u(t) \textnormal{ in } H^1(\mathbb{R}^N)$,
where $u(\cdot;\bar v)$ and $\widehat u$ are as above and the convergence is uniform with respect to $t$ from compact subsets of $[0,+\infty)$.
\end{enumerate}
\end{corollary}

\end{document}